\documentstyle[11pt,twoside]{article}
\newtheorem{theorem}{\noindent\bf Theorem}[section]
\newtheorem{proposition}[theorem]{\noindent\bf Proposition}
\newtheorem{lemma}[theorem]{\noindent\bf Lemma}

\newtheorem{example}[theorem]{\noindent\bf Example}
\newtheorem{remark}[theorem]{\noindent\bf Remark}

\newtheorem{definition}[theorem]{\noindent\bf Definition}

\newcommand{\h}{{\cal H}}
\newcommand{\w}{{\cal W}}
\input{amssym}
\begin{document}
\date{}
\title{{\Large\bf Discretization of continuous frame}}

\author{{\normalsize\sc A. Fattahi, H. Javanshiri}}

\maketitle
\normalsize

\begin{abstract}

In this paper we consider on the  notion of continuous frame of
subspace and define a new concept of continuous frame, entitled
{\it continuous atomic resolution of identity}, for arbitrary
Hilbert space $\h$ which has a countable reconstruction formula.
Among the other result, we characterize the relationship between
this new concept and other known continuous frame. Finally, we
state and prove the assertions of the stability of perturbation in
this concept. \footnote{2000 {\it Mathematics Subject
Classification}: 42C15, 46C99, 94A12, 46B25, 47A05.

   {\it Key words}: Bonded operator, Hilbert space, continuous frame, atomic resolution of identity.}
\end{abstract}

\section{\large\bf Introduction and Preliminaries}

As we know frames are more flexible tools to translate information
than bases, and so they are suitable replacement for bases in a
Hilbert space $\h$.
 Finding a representation of $f\in\h$ as a linear
combination of the elements in frames, is the main goal of
discrete  frame theory. But in continuous frame, which is a
natural generalization from
 discrete, it is
not straightforward.
 However, one
of the applications of frames is in wavelet theory. The practical
implementation of the wavelet transform  in signal processing
requires the selection of a discrete set of points in the
transform space. Indeed, all formulas must generally be evaluated
numerically, and a computer is an intrinsically discrete object.
But this operation must be performed in such a way that no
information is lost. So efforts have been done to find methods to
discretize classical continuous frames for use in applications
like signal processing, numerical solution of PDE, simulation, and
modelling; see for example [1, 8]. In particular, the discrete
wavelet transform and Gabor frames are prominent examples and have
been proven to be a very successful tool for certain applications.
Since the problem of discretization is so important it would be
nice to have a general method for this purpose. For example, Ali,
Antoine, and Gazeau in [1] asked for conditions which ensure that
a certain sampling of a continuous frame $\{\psi_x\}_{x\in X}$
yields a discrete frame $\{\psi_{x_i}\}_{i\in I}$ (see also [9]).

In the recent years, there has been shown considerable interest by
harmonic and functional analysts in the frame of subspace problem
of the separable Hilbert space; see \cite{ck}, \cite{ak}, \cite{a}
and \cite{af} and the references there. Frame of subspace was
first introduced by P. Casazza and G. Kutyniok in \cite{ck}. They
present a reconstruction formula $f=\sum_{i\in
I}\nu_i^2S^{-1}\pi_{W_i}(f)$ for frames of subspace. Continuous
frame of subspace is a natural generalization from
 discrete frame of subspace to continuous.

As we expect, in discrete frame of subspace every element in $\h$
has an expansion in terms of frames. But in the continuous case it
respect to Bochner integral which is not desirable. Therefore,
discretization of continuous frame of subspace is also very
important.


Suppose that the measure $\mu$, which appears in the integral of
continuous frame, is Radon or discontinuous (Note that there exist
infinite many positive finite discontinuous measure on a locally
compact space $X$ which are not counting measure). Then $\{x\in X:
\mu(\{x\}\neq 0\}$ is nonempty set and we may investigate about
some conditions under which every fixed element $f\in\h$ has a
countable subfamily $J_f$ of $X$ with frame property for $h$. This
leads us to define {\it uca-resolution of identity} (Definition
2.1),
 which is a
generalization of the resolution of identity (\cite{ck},
Definition 3.24), and atomic resolution of identity (\cite{ak}),
to arbitrary Hilbert space (separable or nonseparable). We then
show that in this concept many basic properties of discrete state
can be derived within this more general context. In fact
uca-resolution identity helps us to investigate continuous frames
which have discretization. Because under some extra conditions,
every uca-resolution of identity provides a continuous frame of
subspace, and conversely. This means that the relationship between
uca-resolution of identity and known continuous frames, such as
frame of subspace is very tight.

Assume $\h$ to be a Hilbert space and $X$ be a locally compact
Hausdorff space endowed with a positive Radon or discontinuous
measure $\mu$. Let $\w=\{W_x\}_{x\in X}$ be a family of closed
subspaces in $\h$ and let $\omega:X\rightarrow [0,\infty)$ be a
measurable mapping such that $\omega\neq 0$ almost everywhere
(a.e.). We say that $\w_\omega=\{(W_x,\omega(x))\}_{x\in X}$ is a
continuous frame of subspace for
$\h$, if;\\\\
(a) the mapping $x\mapsto \pi_{W_x}$ is weakly measurable;\\
(b) there exist constants $0<A,B<\infty$ such that
$$A\|f\|^2\leq\int_X{\omega(x)}^2\|\pi_{W_x}(f)\|^2\;d\mu(x)\leq
B\|f\|^2~~~~~(1)$$
 for all $f\in\h$.
 The numbers $A$ and $B$ are
called the continuous frame of subspace bounds. If $\w_\omega$
satisfies only the upper inequality in $(1)$, then we say that it
is a continuous Bessel frame of subspace with bound $B$. Note that
if $X$ is a countable set and $\mu$ is the counting measure, then
we obtain the usual definition of a (discrete) frame of subspace.

For each continuous Bessel frame of subspace
$\w_\omega=\{(W_x,\omega(x))\}_{x\in X}$, if we define the
representation space associated with $\w_\omega$ by
$L^2(X,\h,\w_\omega)=\{\varphi:X\rightarrow \h|~\varphi
~\hbox{is~measurable},\\~\varphi (x)\in W_x~\hbox{and}~\int_X
\|\varphi(x)\|^2\;d\mu(x)<\infty\}$, then $L^2(X,\h,\w_\omega)$
with the inner product given by
$$\big<\varphi,\psi\big>=\int_X \big<\varphi(x),\psi(x)\big>\;d\mu(x),~~~~~~~\hbox{for all}~ \varphi, \psi\in L^2(X,\h,\w_\omega),$$
is a Hilbert space. Also, the synthesis operator
$T_{\w_\omega}:L^2(X,\h,\w_\omega)\rightarrow \h$ is define by
$$\big<T_{\w_\omega}(\varphi),f\big>=\int_X\omega(x)\big<\varphi(x),f\big>\;d\mu(x),$$
for all $\varphi\in L^2(X,\h,\w_\omega)$ and $f\in\h$. Its adjoint
operator is $T_{\w_\omega}^*:\h\rightarrow L^2(X,\h,\w_\omega)$;
$T_{\w_\omega}^*(f)=\omega\pi_{\w_\omega}(f)$. For more details
see \cite{af}.

Now, we give two immediate consequences from the above discussion.
As the first, we have the following characterization of continuous
Bessel frame of subspace in term of their synthesis operators as
in discrete frame theory; see \cite{a}.

\begin{theorem}
A family $\w_\omega$ is a continuous Bessel frame of subspace
with Bessel fusion bound $B$ for $\h$ if and only if the synthesis
operator $T_{\w_\omega}$ is a well-defined bounded operator and
$\|T_{\w_\omega}\|\leq\sqrt{B}$.
\end{theorem}

Also, by an argument similar to the proof of (\cite{a}, Theorem
2.6), we have a characterization of continuous frame of subspace
as follows;

\begin{theorem}
The following
 conditions are equivalent:\\\\
{\rm(a)} $\w_\omega=(\{W_x\}_{x\in X},\omega(x))$ is a continuous
frame of subspace for $\h$;\\
{\rm(b)} The synthesis operator $T_{\w_\omega}$ is a bounded
,linear operator from $L^2(X,\h,\w_\omega)$ onto $\h$;\\
{\rm(c)} The analysis operator $T^*_{\w_\omega}$ is injective with
closed range.
 \end{theorem}

If $\w_\omega$ is a continuous frame of subspace for $\h$ with
frame bounds $A, B$, then we define the frame of subspace operator
$S_{\w_\omega}$ for $\w_\omega$ by
$$S_{\w_\omega}(f)=T_{\w_\omega}T^*_{\w_\omega}(f),~~~~~~~~f\in\h,$$
which is a positive, self-adjoint, invertible operator on $\h$ with $A\cdot Id_\h\leq S_{\w_\omega}\leq B\cdot Id_\h.$

\section{\large\bf Main result}

For instituting a relationship between discrete and continuous
frame of subspace, we generalize the concept of continuous frame
and resolution of identity to arbitrary Hilbert space $\h$. For
this propose, we introduce the summation to noncountable form. Let
$\h$ be a Hilbert space and $\{T_x\}_{x\in X}$ be a family of
bounded operators on it. If now, set $\Gamma$, the collection of
all finite subset of $X$, then $\Gamma$ is a directed set ordered
under inclusion.

Let $f$ be a fixed element of the Hilbert space $\h$. Define the
sum $S(f)$ of the family $\{T_x(f)\}_{x\in X}$ as the limit
$$S(f)=\sum_{x\in X}T_x(f)=\lim \{\sum_{x\in \gamma}T_x(f):~\gamma\in\Gamma\}.$$
If this limit exists, we say that the family $\{T_x(f)\}_{x\in X}$
is unconditionally summable. It is easy to see that the family
$\{T_x(f)\}_{x\in X}$ is unconditionally summable if and only if
for each $\varepsilon>0$, there exist a finite subset
$\gamma_0\in\Gamma$ such that
$$\|\sum_{x\in \gamma_1}T_x(f)-\sum_{x\in \gamma_2}T_x(f)\|<\varepsilon,$$
for each $\gamma_1, \gamma_2>\gamma_0$. Therefore for each
$\varepsilon>0$, there is a finite subset $\gamma_0$ of $X$ such
that
$$\|T_x(f)\|<\varepsilon$$
for all $x\in X\setminus\gamma_0$. Hence for a fixed element
$f\in\h$, if $\{T_x(f)\}_{x\in X}$ is unconditionally summable,
then
  $J_f = \{x\in
X:~T_x(f)\neq 0\}$ is countable.

\begin{definition}
Let $\h$ be a Hilbert space and let $\omega:X\rightarrow
[0,\infty)$ be a measurable mapping such that $\omega\neq 0$
almost everywhere. We say that a family of bounded operator
$\{T_x\}_{x\in X}$ on $\h$ is an unconditional continuous atomic
resolution {\rm(}uca-resolution{\rm)} of the identity with respect
to $\omega$ for $\h$, if there exist positive real numbers $C$ and
$D$ such that
for all $f\in\h$,\\\\
{\rm(a)} ~the mapping $x\mapsto T_{x}$ is weakly measurable;\\
{\rm(b)}~$ C\|f\|^2\leq\int_{X}\omega(x)^2\| T_x(f)\|^2d\mu(x)\leq
D\| f\|^2;$\\
{\rm(c)}~$f=\sum_{x\in X} T_x(f).$
\end{definition}
The optimal values of C and D are called the uca-resolution of the
identity bounds. It follows from the definition and the uniform
boundedness principle
that ${\rm sup}_{x\in X}\|T_x\|_{x\in X}\\<\infty$.\\

\begin{remark}
  {\rm(a) If $f\in\h$ satisfies in (c), then as we mention in above,
   there is a countable measurable subset $J_f$ (depends of $f$) of $X$ such that
$$T_x(f)=0,$$
for all $x\in X\setminus J_f.$ So $$\int_X \omega(x)^2\|
T_x(f)\|^2d\mu(x)=\sum_{j\in J_f}\omega(j)^2\|
T_j(f)\|^2\;\mu(\{j\})$$ and condition (b) transform to

$$ C\|f\|^2\leq\sum_{j\in J_f}\omega(j)^2\|
T_j(f)\|^2\;\mu(\{j\})\leq D\| f\|^2$$ (b) If $\h$ is a separable
Hilbert space with an orthonormal bases $\{e_n\}_{n=1}^{\infty}$,
then by condition (c), for each $n$ there exists a countable
measurable subset $J_n$ of $X$ such that
$$T_x(e_n)=0,$$
for all $x\in X\setminus J_n.$ So, we can find a countable subset
$J=\bigcup_{n=1}^\infty J_n$ of $X$ such that
$$T_x(f)=0,$$
for all $f\in\h$ and $x\in X\setminus J$, and we have
$$\int_X \omega(x)^2\| T_x(f)\|^2d\mu(x)=\sum_{j\in J}\omega(j)^2\| T_j(f)\|^2\;\mu(\{j\}).$$
Therefore, if $\h$ is a separable Hilbert space,  Definition 2.1
and  Definition 3.1 in \cite{ak} coincide. }
\end{remark}


From now on $\h$ is a Hilbert space with orthonormal bases
$\{e_\lambda\}_{\lambda\in\Lambda}$ and $X$ is a locally compact
Hausdorff space endowed with a positive Radon or discontinuous
measure $\mu$, and $\omega:X\rightarrow [0,\infty)$ is a
measurable mapping such that $\omega\neq 0$ almost everywhere. For
a fix element $f\in\h$, by [7] there exists a countable subset $J$
of $\Lambda$ such that $\big<f,e_\lambda\big>=0$ for all
$\lambda\in\Lambda\setminus J$.


The following is an important example of uca-resolution compatible
with definition 2.1, and note that this example does not satisfy
in the definition of resolution of identity and atomic resolution
of identity which stated in \cite{ck} and \cite{ak}, respectively.

\begin{example}
{\rm Let $\h$ be a Hilbert space with an orthonormal basis
$\{e_\lambda\}_{\lambda\in\Lambda}$. If, we consider $\Lambda$ as
a locally compact space with discrete topology and measurable
space endowed with counting measure, then the family
$\{T_\lambda\}_{\lambda\in\Lambda}$ of bounded operators on $\h$,
defined by
$$T_\lambda(f)=\big<e_\lambda,f\big>e_\lambda,~~~~~~~~~\hbox{for all}~ f\in\h ~\hbox{and}~ \lambda\in\Lambda,$$
is an uca-resolution of identity for $\h$}.
\end{example}


In the next theorem we show that every uca-resolution of identity
for $\h$, provides a continuous frame of subspace.

\begin{theorem}
Let $\{T_x\}_{x\in X}$ be a family of bounded operators on $\h$
and for each $x\in X$, set $W_x=\overline{T_x(\h)}$. Suppose that
there exists $D>0$ and $R>0$ such that the following conditions
holds:\\
{\rm(a)} $f=\sum_{x\in X}\omega(x)^2T_x(f)\mu(\{x\});$\\
{\rm(b)} $\int_X\omega(x)^2\|\pi_{W_x}(f)-T_x(f)\|^2d\;\mu(x)\leq
R\|f\|^2;$\\
{\rm(c)} $\int_X\omega(x)^2\|T_x(f)\|^2d\;\mu(x)\leq D\|f\|^2,$\\
for all $f\in\h$. Then $\{(W_x,\omega(x))\}_{x\in X}$ is a
continuous frame of subspace for $\h$.
\end{theorem}
{\noindent Proof.} Let $f$ be a fix element of $\h$. as we mention
in remark 2.2(a), there exists a countable subset $J_f$ of $X$
such that
$$\omega(x)^2T_x(f)\mu(\{x\})=0,$$
for all $x\in X\setminus J_f$, and
$$\int_X\omega(x)^2\|T_x(f)\|^2\;d\mu(x)=\sum_{x\in X}\omega(x)^2\|T_x(f)\|^2\mu(\{x\}).$$
So we can use Cauchy-Schwarz inequality and compute as follows
\begin{eqnarray*}
\|f\|^4&=&(\langle f,\sum_{x\in X}\omega(x)^2T_x(f)\mu(\{x\})\rangle)^2\\
&=&(\sum_{x\in
X}\omega(x)\langle\sqrt{\mu(\{x\})}\;f,\omega(x)\sqrt{\mu(\{x\})}\;T_x(f)\rangle)^2\\
&=&(\sum_{x\in
X}\omega(x)\langle\sqrt{\mu(\{x\})}\;\pi_{W_x}(f),\omega(x)\sqrt{\mu(\{x\})}\;T_x(f)\rangle)^2\\
&\leq&(\sum_{x\in
X}\omega(x)\|\sqrt{\mu(\{x\})}\;\pi_{W_x}(f)\|\|\omega(x)\sqrt{\mu(\{x\})}\;T_x(f)\|)^2\\
&\leq&(\sum_{x\in
X}\omega(x)^2\|\pi_{W_x}(f)\|^2\mu(\{x\}))(\sum_{x\in X}\|\omega(x)\sqrt{\mu(\{x\})}\;T_x(f)\|^2)\\
&\leq&(\int_{x\in
X}\omega(x)^2\|\pi_{W_x}(f)\|^2\;d\mu(x))(\int_X\omega(x)^2\|T_x(f)\|^2\;d\mu(x))\\
&\leq&D\|f\|^2(\int_X\omega(x)^2\|\pi_{W_x}(f)\|^2\;d\mu(x)).
\end{eqnarray*}
Also, by triangle inequality and hypothesis we have
\begin{eqnarray*}
\int_X\omega(x)^2\|\pi_{W_x}(f)\|^2\;d\mu(x)\leq
D(1+\sqrt{\frac{R}{D}})^2\|f\|^2,
\end{eqnarray*}
so the assertion holds.$\hfill\Box$\\


Casazza and Kutyniok in \cite{ck} introduced an interesting
example of atomic resolution of identity. In the next theorem we
obtain the uca-resolution of identity form, which is the converse
of theorem 2.4.

\begin{theorem}
Let $\{(W_x,\omega(x))\}_{x\in X}$ be a continuous Bessel frame of
subspace for $\h$ with Bessel bound $D$, and for each $x\in X$,
let $T_x:\h\rightarrow W_x$ be a bounded operator such that
$T_x\pi_{W_x}=T_x$. Also assume that for each $f\in\h$
$$f=\sum_{x\in X}\omega(x)^2T_x(f)\mu(\{x\}).$$
Then for all $f\in\h$ we have
$$\frac{1}{D}\|f\|^2\leq\int_X\omega(x)^2\|T_x(f)\|^2\;d\mu(x)\leq DE\|f\|^2,$$
where $E={\rm sup}_{x\in X}\|T_x\|_{x\in X}$.
\end{theorem}
{\noindent Proof.} By similar prove of Theorem 2.4, we obtain
$$\frac{1}{D}\|f\|^2\leq\int_X\omega(x)^2\|T_x(f)\|^2\;d\mu(x).$$
Also we have
\begin{eqnarray*}
\frac{1}{D}\|f\|^2&\leq&\int_X\omega(x)^2\|T_x(f)\|^2\;d\mu(x)\\
&=&\int_X\omega(x)^2\|T_x\pi_{W_x}(f)\|^2\;d\mu(x)\\
&\leq&\int_X\omega(x)^2\|T_x\|^2\|\pi_{W_x}(f)\|^2\;d\mu(x)\\
&\leq&E\int_X\omega(x)^2\|\pi_{W_x}(f)\|^2\;d\mu(x)\leq DE\|f\|^2.
\end{eqnarray*}
Whence, for each $f\in\h$
$$\frac{1}{D}\|f\|^2\leq\int_X\omega(x)^2\|T_x(f)\|^2\;d\mu(x)\leq DE\|f\|^2.$$
as we required.$\hfill\Box$\\

~~~~~~~~~~~~~~~~~~~~~~~~~~~~~~~~~~~~~~~~~~~~~~~~~~~~~~~~~~~~~~~~~~~~~~~~~~~~~~~~~~~~~~~~~~~~~~~~~~~~~~~~~~~~~~~~~~~~~~~~~~~~~~~~~~~~~~~~~~~~~

\begin{proposition}
Let $\{W_x\}_{x\in X}$ be a family of closed subspace of Hilbert
space $\h$ such that the mapping $x\rightarrow \pi_{W_x}$ is
weakly
measurable. Also suppose $\omega$ is a bounded map and the following conditions hold for all $f\in\h$:\\
{\rm(a)} There exists $C>0$ such that
$$\int_X\parallel\pi_{W_x}(f)\parallel^2d\mu(x)\leq \frac{1}{C}\parallel
f\parallel^2,$$\\
{\rm(b)} $f=\sum_{x\in X}\omega(x)\pi_{W_x}(f)\mu(\{x\})$.\\
Then $\{(W_x,\omega(x))\}_{x\in X}$ is a continuous frame of
subspace for $\h$.
\end{proposition}

{\noindent{Proof.}} By condition (a) we see that
$$\int_X\omega(x)^2\|\pi_{W_x}(f)\|^2\;d\mu(x)\leq\frac{\sup_{x\in X}\omega(x)}{C}\|f\|^2, ~~~~~~~(f\in\h).$$
Condition (b) implies that for a fixed element $f$ of $\h$
$$\int_X\omega(x)^2\|\pi_{W_x}(f)\|^2\;d\mu(x)=\sum_{x\in X}\omega(x)^2\|\pi_{W_x}(f)\|^2\mu(\{x\}),$$
and
$$\int_X\|\pi_{W_x}(f)\|^2\;d\mu(x)=\sum_{x\in X}\|\pi_{W_x}(f)\|^2\mu(\{x\}).$$
Now, since the family $\{\omega(x)\mu(\{x\})T_x\}$ is
unconditional summable, we can use Cauchy-Schwarz inequality and
compute as follows

\begin{eqnarray*}
\|f\|^4&=&(\langle\sum_{x\in
X}\omega(x)\mu(\{x\})\pi_{W_x}(f),f\big>)^2\\
&=&(\sum_{x\in X}\omega(x)\mu(\{x\})\|\pi_{W_x}(f)\|^2)^2\\
&\leq&(\sum_{x\in
X}\omega(x)^2\mu(\{x\})\|\pi_{W_x}(f)\|^2)(\sum_{x\in
X}\|\pi_{W_x}(f)\|^2\mu(\{x\}))\\
&\leq& \frac{1}{C}\|f\|^2(\sum_{x\in
X}\omega(x)^2\|\pi_{W_x}(f)\|^2\mu(\{x\}))
\end{eqnarray*}
Thus
$$C\|f\|^2\leq\sum_{x\in X}\omega(x)^2\|\pi_{W_x}(f)\|^2\mu(\{x\})=\int_X\omega(x)^2\|\pi_{W_x}(f)\|^2\;d\mu(x)$$
for all $f\in\h$, and this complete the proof.$\hfill\Box$\\


In the following proposition we give a reconstruction formula for
continuous frame of subspace in the special case.

\begin{proposition}
Let $\{W_x\}_{x\in X}$ be a family of orthogonal closed subspace
of Hilbert space $\h$. If $\{(W_x,\omega(x)\}_{x\in X}$ is a
continuous frame of subspace for $\h$ with bounds $C, D$, then for
each $f\in\h$
$$f=\sum_{x\in X}\pi_{W_x}(f).$$
The converse is true if $\omega$ is bounded and there exists $C>0$
such that
$$\int_X\|\pi_{W_x}(f)\|^2\;d\mu(x)\leq\frac{1}{C}\|f\|^2,$$
for all $f\in\h$.
\end{proposition}

{\noindent Proof.} Let $\{(W_x,\omega(x)\}_{x\in X}$ be a
continuous frame of subspace. First, we should noted that for each
$f\in\h$, by Hahn-Banach Theorem and orthogonality of the family
$\{W_x\}_{x\in X}$, there exists a sequence $\{f_n\}$ in $\h$ such
that $f_n\longrightarrow f$ and for each $n$ we have the following
equality
$$f_n=\sum_{x\in X}\pi_{W_x}(f_n).$$

Now we define $S_\gamma(f)=\sum_{x\in\gamma}\pi_{W_x}(f)$, where
$\gamma$ is an arbitrary finite subset of $X$ and $f\in\h$.
Therefore
\begin{eqnarray*}
C\|S_\gamma(f)-f\|^2&\leq&\int_X\omega(x)^2\|\pi_{W_x}(S_\gamma(f)-f)\|^2\;d\mu(x)\\
&\leq&\int_X\omega(x)^2\|\pi_{W_x}(f)\|^2\;d\mu(x)\\
&\leq&D\|f\|^2.
\end{eqnarray*}
By replacing $f$ with $f_n-f$ we obtain
$$\|S_\gamma(f_n-f)-(f_n-f)\|\leq\sqrt{\frac{D}{C}}\|f_n-f\|.$$
 The converse holds by 2.6.$\hfill\Box$\\


Now we want to show that, by a given uca-resolution of identity,
each $f\in\h$ has a new countable reconstruction formula. First we
need the following Lemma:
\begin{lemma}
Let $\{T_x\}_{x\in X}$ be an uca-resolution of the identity with
respect to weight $\omega$ for $\h$ with bounds $C$ and $D$, and
let $\{f_i\}_{i\in I}$ be a frame sequence. Then there exists a
countable subset $J$ of $X$, such that $\{\omega(j)\sqrt{\mu
(\{j\})}T^*_j(f_i)\}_{i\in I, j\in J}$ is a frame for
$\overline{\hbox{span}}\{f_i\}_{i\in I}$.
\end{lemma}
{\noindent {Proof.}} If we set $J_i=\{x\in X:~T_x(f_i)\neq 0\}$,
then by definition of uca-resolution of identity , $J_i$ is a
countable and measurable subset of $X$. Now, set $J=\bigcup_{i\in
I}J_i$. So $J$ is a countable and measurable subset of $X$, and
for each $f\in\overline{\hbox{span}}\{f_i\}_{i\in I}$ and $x\in
X\setminus J$ we have
$$T_x(f)=0.$$
Hence we see that for each
$f\in\overline{\hbox{span}}\{f_i\}_{i\in I}$
$$C\|f\|^2\leq \sum_{j\in J}\omega^2(j)\mu(\{j\})\|T_j(f)\|^2\leq D\|f\|^2,$$
and
$$f=\sum_{j\in J} T_j(f),$$
and these series converge unconditionally.

Now, suppose that $A$ and $B$ are frame bounds of $\{f_i\}_{i\in
I}$. For each $f\in\overline{\hbox{span}}\{f_i\}_{i\in I}$ we have
$$A\sum_{j\in J}\omega^2(j)\mu(\{j\})\|T_j(f)\|^2\leq\sum_{j\in J}\sum_{i\in I}|<\omega^2(j)\mu(\{j\})T_j(f),f_i>|^2$$
$$~~~~~~~~~~~~~~~~~~~~~~~\leq B \sum_{j\in J}\omega^2(j)\mu(\{j\})\|T_j(f)\|^2,$$
and therefore
$$AC\|f\|^2\leq A\sum_{j\in J}\omega^2(j)\mu(\{j\})\|T_j(f)\|^2$$
$$~~~~~~~~~~~~~~~~~~~~~~~~\leq \sum_{j\in J}\sum_{i\in
I}|<f,\omega^2(j)\mu(\{j\})T^*_j(f_i)>|^2$$
$$~~~~~~~~~~~~~~~~~~~~~~~~~~~\leq B\sum_{j\in J}\omega^2(j)\mu(\{j\})\|T_j(f)\|^2\leq BD\|f\|^2.$$
and this complete the proof.$\hfill\Box$\\


\begin{theorem}
Let $\{T_x\}_{x\in X}$ be an uca-resolution of the identity with
respect to weight $\omega$ for $\h$ with bounds $C$ and $D$. Then
for each $f\in\h$, there exists a countable subset $I$
{\rm(}dependents on $f${\rm)} of $X$, such that we have the
following reconstruction formula
$$f=\sum_{i\in I}\omega^2(i)\mu(\{i\})S^{-1}T^*_iT_i(f)=\sum_{i\in I}\omega^2(i)\mu(\{i\})T^*_iT_iS^{-1}(f),$$
where $S$ is a frame operator of a frame sequence.
\end{theorem}
{\noindent Proof.} Let $f$ be a fix element of Hilbert space $\h$.
Set
$$\h_f=\overline{\hbox{span}}\{e_j\}_{j\in J},$$
where $J=\{j\in\Lambda:~\big<e_j,f\big>\neq 0\}$ is a countable
subset of $\Lambda$. Then, by Lemma 2.8, there is a countable
subset $I$ of $X$ such that the sequence $\{\omega(i)\sqrt{\mu
(\{i\})}T^*_i(e_j)\}_{i\in I, j\in J}$ is a
frame for $\h_f$.\\
If now, $S\in B(\h)$ is the frame operator of
$\{\omega(i)\sqrt{\mu (\{i\})}T^*_i(e_j)\}_{i\in I, j\in J}$, then
we have
$$S(f)=\sum_{i\in I}\sum_{j\in J}\big<f,\omega(i)\sqrt{\mu (\{i\})}T^*_i(e_j)\big>\omega(i)\sqrt{\mu (\{i\})}T^*_i(e_j)$$

$$=\sum_{i\in I}\omega^2(i)\mu(\{i\})T^*_i(\sum_{j\in J}\big<T_i(f),e_j\big>e_j)~~~~~~~~~~~$$
$$=\sum_{i\in I}\omega^2(i)\mu(\{i\})T^*_iT_i(f).~~~~~~~~~~~~~~~~$$
Hence, the reconstruction formula follows immediately from the invertibility of the operator $S$.$\hfill\Box$\\


In the rest of paper we consider to stability of perturbation in
uca-resolution of identity. First, let us state and proof of the
following useful lemma.

\begin{lemma}
Let $\{T_x\}_{x\in X}$ and $\{S_x\}_{x\in X}$ be two families of
bounded operators on $\h$ and there exists $0<\lambda<1$ such that
for all finite subset $I$ of $X$
$$\|\sum_{i\in I}(T_i-S_i)(f)\|\leq\lambda\|\sum_{i\in I}T_i(f)\|~~~~~~{\rm(}f\in\h{\rm)}~~~~{\rm(}1{\rm)}.$$
If $\{(T_x,\omega(x)\}_{x\in X}$ is an uca-resolution of identity
then we have the following reconstruction formula
$$f=\sum_{x\in X}S_xS^{-1}(f)~~~~~~{\rm(}f\in\h{\rm)}$$
where $S$ is an invertible operator on $\h$.
\end{lemma}
{\noindent Proof.} Let $f\in\h$ and let $I$ be a finite subset of
$X$. Since
$$\|f-\sum_{i\in I}S_i(f)\|\leq\|f-\sum_{i\in I}T_i(f)\|+\|\sum_{i\in I}T_i(f)-\sum_{i\in I}S_i(f)\|.$$
Therefore by inequality (1) we have
$$\|f-\sum_{i\in I}S_i(f)\|\leq\|f-\sum_{i\in I}T_i(f)\|+\lambda\|\sum_{i\in I}T_i(f)\|~~~~~~{\rm(}2{\rm)}.$$
Hence, the family $\{S_x(f)\}_{x\in X}$ is unconditionally
summable. Now, we define $S:\h\rightarrow\h$ by $S(f)=\sum_{x\in
X}S_x(f)$. By inequality (2) and using that $\{(T_x,\omega(x)\}$
is assumed to be uca-resolution of identity,  $S$ is well defined
and we have
$$\|f-S(f)\|\leq\lambda\|f\|,$$
for all $f\in\h$. So $\|{\rm id}_\h-S\|\leq\lambda<1$, and
therefore $S$ is an invertible operator on $\h$. Hence for all
$f\in\h$ we have
$$\sum_{x\in X}S_xS^{-1}(f)=SS^{-1}(f)=f,$$
and this complete the proof.$\hfill\Box$\\


\begin{definition}
Let $\{T_x\}_{x\in X}$ and $\{S_x\}_{x\in X}$ be two families of
bounded operators on  $\h$, and let $\omega:X\rightarrow
[0,\infty)$ be measurable map such that $\omega(x)\neq 0$ almost
everywhere. Suppose that  $0\leq\lambda_1, \lambda_2<1$, and
$\varphi:X\rightarrow [0,\infty)$ is an arbitrary positive map
such that $\int_X\varphi(x)^2\;d\mu(x)<\infty$. If
$$\|\omega(x)(T_x-S_x)(f)\|\leq\lambda_1\|\omega(x)T_x(f)\|+\lambda_2\|\omega(x)S_x(f)\|+\varphi(x)\|f\|$$
for all $f\in\h$ and $x\in X$, then we say that
$\{(S_x,\omega(x))\}_{x\in X}$ is a $(\lambda_1, \lambda_2,
\varphi)$-perturbation of $\{(T_x,\omega(x))\}_{x\in X}$.
\end{definition}


From now on let $\{S_x\}_{x\in X}$ be a family of bounded
operators on $\h$ such that the mapping $x\mapsto S_{x}(f)$ is
weakly measurable. Then for each bounded operator
$S:\h\rightarrow\h$, the map $x\mapsto S_{x}S(f)$ is weakly
measurable. Hence by Lemma 2.9, we have the following theorem.

\begin{theorem}
Let $\{(T_x,\omega(x))\}_{x\in X}$ be an uca-resolution of
identity for $\h$ with bounds $C$ and $D$, and let
$\{(S_x,\omega(x))\}_{x\in X}$ be a $(\lambda_1, \lambda_2,
\varphi)$-perturbation of $\{(T_x,\omega(x))\}_{x\in X}$ for some
$0\leq\lambda_1,\lambda_2<1$. Moreover assume that
$(1-\lambda_1)\sqrt{C}-(\int_X\varphi(x)^2\;d\mu(x))^{\frac{1}{2}}>0$
and for some $0\leq\lambda<1$
$$\|\sum_{i\in I}(T_i-S_i)(f)\|\leq\lambda\|\sum_{i\in I}T_i(f)\|~~~~~~~~~{\rm(}f\in\h{\rm)},$$
for all finite subset $I$ of $X$. Then there exist an invertible
operator $S$ on $\h$ such that $\{(S_xS^{-1},\omega(x)\}_{x\in X}$
is a uca-resolution of the identity on $\h$.
\end{theorem}
{\noindent Proof.} First it should be noted that by Lemma 2.10,
there exists an invertible operator $S$ on $\h$, such that the
family $\{S_xS^{-1}\}_{x\in X}$ satisfies in 2.1(c). Also by Open
mapping Theorem and Closed Graph Theorem, there exist $A>0$ and
$B>0$ such that
$$A\|f\|\leq\|S^{-1}(f)\|\leq B\|f\|$$
for all $f\in\h$.\\

 Now, for  $f\in\h$ we obtain\\\\
$(\int_X\omega(x)^2\|S_x(f)\|^2\;d\mu(x))^\frac{1}{2}\leq(\int_X\omega(x)^2(\|T_x(f)\|+\|(T_x-S_x)(f)\|)^2\;d\mu(x))^\frac{1}{2}$\\\\
$\leq(\int_X((\omega(x)^2(\|T_x(f)\|+\lambda_1\|T_x(f)\|\|+\lambda_2\|S_x(f)\|))+\varphi(x)\|f\|)^2\;d\mu(x))^\frac{1}{2}$\\\\
$\leq(1+\lambda_1)(\int_X\omega(x)^2\|T_x(f)\|^2\;d\mu(x))^\frac{1}{2}
+\lambda_2(\int_X\omega(x)^2\|S_x(f)\|^2\;d\mu(x))^\frac{1}{2}
+\|f\|(\int_X\varphi(x)^2\;d\mu(x))^\frac{1}{2}.$\\\\
 Hence
 $$\int_X\omega(x)^2\|S_xS^{-1}(f)\|^2\;d\mu(x)\leq(\frac{(1+\lambda_1)\sqrt{D}+(\int_x\varphi(x)^2\;d\mu(x))^\frac{1}{2}}{1-\lambda_2})^2B^2\|f\|^2.$$\\\\
To prove the lower bound, first we observe that
$$\|f\|^2\leq\frac{1}{C}\int_X\omega(x)^2\|T_x(f)\|^2\;d\mu(x),$$
for all $f\in\h$. Therefore, by triangle inequality we have\\\\
$(\int_X\omega(x)^2\|T_x(f)\|^2\;d\mu(x))^\frac{1}{2}-(\int_X\omega(x)^2\|S_x(f)\|^2\;d\mu(x))^\frac{1}{2}
\leq(\int_X\|\omega(x)(T_x-S_x)(f)\|^2)^\frac{1}{2}$\\\\
$\leq\lambda_1(\int_X\omega(x)^2\|T_x(f)\|^2\;d\mu(x))^\frac{1}{2}
+\lambda_2(\int_X\omega(x)^2\|S_x(f)\|^2\;d\mu(x))^\frac{1}{2}
\\\\+\frac{1}{\sqrt{C}}(\int_X\varphi(x)^2\;d\mu(x))^\frac{1}{2}(\int_X\omega(x)^2\|T_x(f)\|^2\;d\mu(x))^\frac{1}{2}.$\\\\
Hence
$$(\frac{1-\lambda_1-\frac{1}{\sqrt{C}}(\int_x\varphi(x)^2\;d\mu(x))^\frac{1}{2}}{1+\lambda_2})
(\int_X\omega(x)^2\|T_x(f)\|^2\;d\mu(x))^\frac{1}{2}\leq(\int_X\omega(x)^2\|S_x(f)\|^2\;d\mu(x))^\frac{1}{2}.$$
So
$$(\frac{(1-\lambda_1)\sqrt{C}-(\int_x\varphi(x)^2\;d\mu(x))^\frac{1}{2}}{1+\lambda_2})^2A^2\|f\|^2\leq\int_X\omega(x)^2\|S_xS^{-1}(f)\|^2\;d\mu(x),$$\\
as we required.$\hfill\Box$\\


\begin{remark}
{\rm Suppose $\{T_x\}_{x\in X}$ and $\{S_x\}_{x\in X}$ are two
families  of bounded operators on $\h$. If
$\{(T_x,\omega(x))\}_{x\in X}$ is a uca-resolution of identity,
then by Cauchy-Schwarz inequality we have
\begin{eqnarray*}
|\big<T_xS_x(f),g\big>|&=&|\big<S_x(f),T^*_x(g)\big>|\\
&\leq&\|S_x(f)\|\|T^*_x\|\|g\|\\
&\leq&\|S_x(f)\|\|g\|\sup_{x\in X}\|T_x\|,
\end{eqnarray*}
for all $f,g\in\h$ and $x\in X$. Hence, for each $f\in\h$ and
$x\in X$
$$\|T_xS_x(f)\|\leq\|S_x(f)\|E,$$
where $E=\sup_{x\in X}\|T_x\|$. }
\end{remark}

\begin{theorem}
Let $\{(T_x,\omega(x))\}_{x\in X}$ be an uca-resolution of
identity for $\h$ with bounds $C$ and $D$, and let $\{S_x\}_{x\in
X}$ be a family of bounded operators on  $\h$ such that for some
$K$
$$\int_X\omega(x)^2\|S_x(f)\|^2\;d\mu(x)\leq D\|f\|^2,$$
for all $f\in\h$. Suppose that  $\varphi:X\rightarrow [0,\infty)$
is a positive map, and  there exist $0<\lambda_1, \lambda_2<1$
such that
$$\|\omega(x)f-\omega(x)^2T_xS_x(f)\|\leq\lambda_1\|\omega(x)T_x(f)\|+\lambda_2\|\omega(x)^2T_xS_x(f)\|+\varphi(x)\|f\|$$
Also
$$\|\sum_{i\in I}(T_i-S_i)(f)\|\leq\lambda\|\sum_{i\in I}T_i(f)\|$$
for all finite subset $I$ of $X$ and for all $f\in\h$, where
$0<\lambda<1$. If $\int_X\varphi(x)\;d\mu(x)<\infty$ and
$0<(\int_X\omega(x)^2d\mu(x))^{\frac{1}{2}}-\lambda_1\sqrt{D}-(\int_X\varphi(x)^2d\mu(x))<\infty$,
then there exists an invertible operator $S$ on $\h$ such that
$\{(S_xS^{-1},\omega(x))\}_{x\in X}$ is an uca-resolution of the
identity on $\h$.
\end{theorem}
{\it Proof.} For $f\in\h$ we have\\\\
$\|f\|(\int_X\omega(x)^2\;d\mu(x))^{\frac{1}{2}}\leq
(\int_X(\|\omega(x)f-\omega(x)^2T_xS_x(f)\|+\|\omega(x)^2T_xS_x(f)\|)^2\;d\mu(x))^{\frac{1}{2}}$\\\\
$\leq
(\int_X\|\omega(x)f-\omega(x)^2T_xS_x(f)\|^2d\mu(x))^{\frac{1}{2}}+(\int_X\|\omega(x)^2T_xS_x(f)\|^2d\mu(x))^{\frac{1}{2}}$\\\\
$\leq
(\int_X(\lambda_1\|\omega(x)T_x(f)\|+\lambda_2\|\omega(x)^2T_xS_x(f)\|+\varphi(x)\|f\|)^2d\mu(x))^{\frac{1}{2}}$\\\\
$+
(\int_X\|\omega(x)^2T_xS_x(f)\|^2d\mu(x))^{\frac{1}{2}}$\\\\
$\leq\lambda_1\sqrt{D}\|f\|+(1+\lambda_2)(\int_X\omega(x)^2\|T_xS_x(f)\|^2d\mu(x))^{\frac{1}{2}}+\|f\|(\int_X\varphi(x)^2d\mu(x))^{\frac{1}{2}}$\\\\
$\leq\lambda_1\sqrt{D}\|f\|+(1+\lambda_2)E(\int_X\omega(x)^2\|S_x(f)\|^2d\mu(x))^{\frac{1}{2}}+\|f\|(\int_X\varphi(x)^2d\mu(x))^{\frac{1}{2}}$\\\\
where $E=\sup_{x\in X}\|T_x\|$. Therefore
$$\|f\|\frac{(\int_X\omega(x)^2d\mu(x))^{\frac{1}{2}}-\lambda_1\sqrt{D}-(\int_X\varphi(x)^2d\mu(x))^{\frac{1}{2}}}{E(1+\sqrt{\lambda_2})}
\leq(\int_X\omega(x)^2\|S_x(f)\|^2d\mu(x))^{\frac{1}{2}}.$$ Now by
Lemma 2.10, and similar to prove of 2.12, the assertion holds.$\hfill\Box$\\



 \vspace{5mm}

 \noindent (Abdolmajid Fattahi) Department of Mathematics, Razi University,
Kermanshah, Iran.\\
E-mail address: majidzr@razi.ac.ir \& abfattahi@yahoo.ca\\\\
(H. Javanshiri) Department of Mathematical Sciences, Isfahan
University of Technology, Isfahan 84156-83111, Iran.\\
E-mail address: h.javanshiri@math.iut.ac.ir \&
hjavanshirigh@yahoo.com


\footnotesize

\begin{thebibliography}{99}

\bibitem{af} {\sc Ali, S. T., Antoine, J. P. and Gazeau, J. P.}:  Coherent States, Wavelets
and their Generalizations, {\it Springer-Verlag} (2000).


\bibitem{af} {\sc Ahmadi, R. and Faroughi, M. H.}: Some
properties of C-fusion frames,
 {\it Turk J Math} {\bf
33} (2009), 1-23.

\bibitem{a} {\sc Asgari, M. S.}: New characterizations of fusion frames (frames of subspaces),
 {\it Proc. Indian Acad. Sci. (Math. Sci.)} {\bf
119} No. 3, (2009), 369–382.



\bibitem{ak} {\sc Asgari, M. S. and Khosravi, A.}: Frames and bases of subspaces in Hilbert spaces,
 {\it J. Math. Anal. Appl.} {\bf
308} (2005), 541-553.

\bibitem{ck} {\sc Casazza, P. G. and Kutyniok, G.}: Frame of subspaces, Wavelets, Frames and Operator
Theorey,
 {\it Contemp. Math} {\bf
345} (1995), 87-113.

\bibitem{c} {\sc Christensen, O.}: Introduction to frames and
Riesz bases, {\it Boston, Birkhauser 2003}.


\bibitem{conway} {\sc Conway, J. B.}: A course in functional analysis, {\it Springer-verlag, New york Inc} 1985.

\bibitem{conway} {\sc Dahlke, S., Fornasier, M. and Raasch, T.} Adaptive frame methods for
elliptic operator equations, Bericht Nr. 2004-3, Fachbereich
Mathematik und Informatik, Philipps-Universit\"{a}t Marburg, 2004.

\bibitem{conway} {\sc Massimo, F. and Holger, R.} Continuous frames, function spaces, and the discretization
problem. {\it J. Fourier Anal. Appl.} {\bf 11} (2005), no. 3,
245–287.






\end{thebibliography}
       \end{document}